\documentclass{article}
\usepackage{graphicx}
\usepackage{amsmath}
\usepackage{listings}
\usepackage{url}
\usepackage{color}
\usepackage{subfigure}
\usepackage{listings}%python code
\usepackage{xcolor}
\usepackage{amsmath}%matrix
\usepackage{color}

\title{Using Crank-Nikolson Scheme to Solve the Korteweg-de Vries (KdV) Equation \\

\footnote{This is the final project report for Math 206C, UCSB}
}
\author{Qiming Wu \\ qimingwu@cs.ucsb.edu}

\begin{document}

\maketitle

\section{Intuitives and Initial Attempts}
The Korteweg-de Vries (KdV) Equation is a well-known mathematical model that describes the behavior of waves on a shallow water surface \footnote{\url{https://en.wikipedia.org/wiki/Korteweg\%E2\%80\%93De_Vries_equation}}. In this project, I will explore 2 approaches (explicit and implicit finite difference method) to solving this equation and compare their performance.

Given the KdV equation
\begin{equation}
u_t - \frac{3}{2}uu_x + u_{xxx} = 0,
\end{equation}\label{equ:kdv}where $u(x,t)$ is the wave amplitude, $x$ is the spatial coordinate, $t$ is time, and subscripts denote partial derivatives. The equation is defined on an interval and the initial data is
\begin{equation}
    u(x,0) = \frac{c}{8}sech^2(\frac{\sqrt{c}}{2}x),
\end{equation}
where
\begin{equation}
    sech(x) = \frac{2}{e^x + e^{-x}}.
\end{equation}

\begin{figure*}
    \centering
    \includegraphics[width=0.6\linewidth]{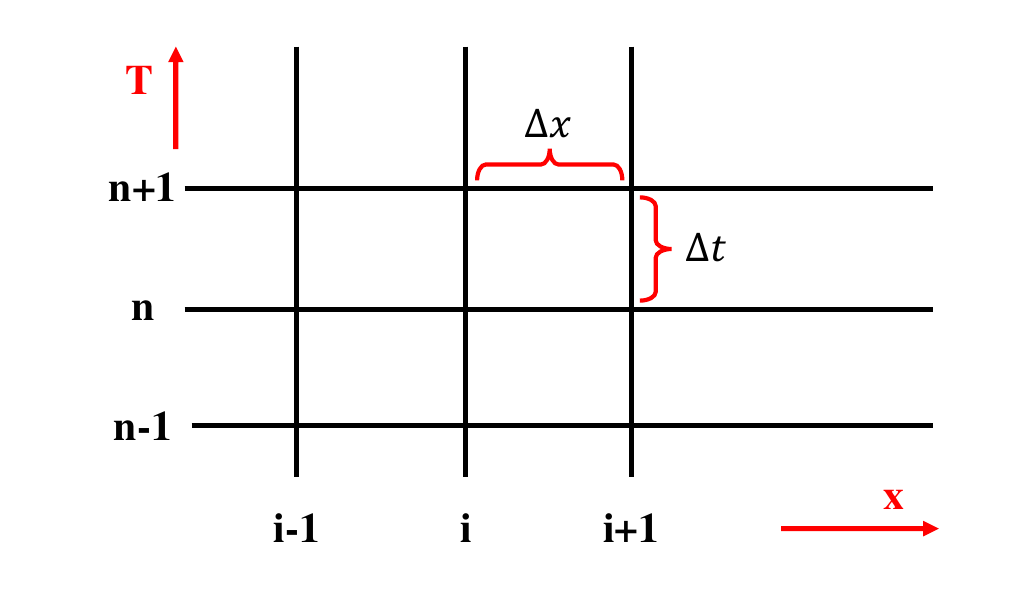}
    \caption{The grid for the partial derivative equations in this paper.}
    \label{fig:grid}
\end{figure*}

Also, to better show the finite difference method, I show the grid for solving this KdV equation in Figure \ref{fig:grid}.

Previously, I only dealt with second-order derivatives in the equation when using the finite difference scheme to solve some partial derivative equations (e.g., the 1D heat equation.). But for this KdV equation, I notice that there are first-order $u_x$ and third-order derivatives $u_{xxx}$ in it and they make things difficult.

In order to solve this dilemma, I first use a centered difference approximation for the first and third derivatives. To discrete the problem in spaces and time, I consider $x_i = i\Delta x$ and $t_n = n\Delta t$ be the discrete spatial and time coordinates, where $\Delta x$ and $\Delta t$ are the spatial and time steps, respectively. Therefore, using a centered difference approximation and an upwind scheme for the nonlinear term in the equation, I can obtain the following results:

\begin{align}
& u_t \approx \frac{u_i^{n+1} - u_i^n}{\Delta t}, \
\\ & u_x \approx \frac{u_{i+1}^n - u_{i-1}^n}{2\Delta x}, \
\\ & u_{xx} = \frac{u^n_{i+1} - 2u^n_{i} + u^n_{i-1}}{(\Delta x)^2},
\
\\ & u_{xxx} \approx \frac{u_{i+2}^n - 2u_{i+1}^n + 2u_{i-1}^n - u_{i-2}^n}{2\Delta x^3}, \
\\ & uu_x \approx u_i^n\frac{u_{i+1}^n - u_{i-1}^n}{2\Delta x}.
\end{align}

Then I replace the partial derivative parts in the KdV equation with the above finite difference approximations. I get the following results:

\begin{equation}
\frac{u_i^{n+1} - u_i^n}{\Delta t} - \frac{3}{2} u_i^n\frac{u_{i+1}^n - u_{i-1}^n}{2\Delta x} + \frac{u_{i+2}^n - 2u_{i+1}^n + 2u_{i-1}^n - u_{i-2}^n}{2\Delta x^3} = 0.
\end{equation}

Then, moving $u^{n+1}_i$ to the left side and moving other terms to the right side, I get the following expressions:

\begin{align}
    u^{n+1}_i = u^{n}_i[1 + \frac{3\Delta t}{4\Delta x}(u^n_{i+1}-u^n_{i-1})] - \frac{\Delta t}{2\Delta x^3}(u^n_{i+2}-2u^n_{i+1}+2u^n_{i-1}-u^n_{i-2}).
\end{align}

Based on the above equation, we can get the results of $u^{n+1}_i$ by iteratively solving for $u$ at each time step, starting from an initial condition $u(x,0)$. However, this method is an explicit method, which means that it is conditionally stable. To better solve this equation, I would change the method with the Crank-Nikolson Scheme.

\section{Crank-Nikolson Scheme}
Since the explicit method of solving the KdV equation is conditionally stable, I try to find an implicit method to better solve the problem. 
To discretize the equation using the Crank-Nicolson scheme, we let $x_i = i\Delta x$ and $t_n = n\Delta t$ be the discrete spatial and time coordinates, where $\Delta x$ and $\Delta t$ are the spatial and time steps, respectively. Using a central difference approximation for the spatial derivative and a weighted average of the forward and backward time steps for the time derivative, we obtain:
\begin{align}
& u_t \approx \frac{u_i^{n+1} - u_i^n}{\Delta t},
\\ & u_x \approx \frac{u^n_{i+1}-u^n_{i-1}+u^{n+1}_{i+1}-u^{n+1}_{i-1}}{4\Delta x},
\\ & u_{xxx} \approx \frac{(u^n_{i+2}-2u^n_{i+1}+2u^n_{i-1}-u^n_{i-2})+(u^{n+1}_{i+2}-2u^{n+1}_{i+1}+2u^{n+1}_{i-1}-u^{n+1}_{i-2})}{4 \Delta x^3}.
\end{align}

However, one of the difficulties is how to linearize the non-linear term in the KdV equation. That is, how to linearize the $u u_{x}$ in the KdV equation. In the view of implicit central difference, we approximate the terms in the following ways:

\begin{align}
& uu_x \approx \textcolor{red}{u_i^n} \frac{u^n_{i+1}-u^n_{i-1}+u^{n+1}_{i+1}-u^{n+1}_{i-1}}{4\Delta x},
\end{align}\label{equ:appro_un}
or
\begin{align}
uu_x \approx \textcolor{red}{u_i^{n+1}} \frac{u^n_{i+1}-u^n_{i-1}+u^{n+1}_{i+1}-u^{n+1}_{i-1}}{4\Delta x}.
\end{align}\label{equ:appro_un+1}
Then, I substitute these approximations into the KdV Equation \ref{equ:kdv} and obtain the following results.

First, I define some terms in the following:
\begin{align}
    \frac{\Delta t}{\Delta x^3} = \alpha,\  \frac{\Delta t}{\Delta x} = \beta
\end{align}

If I use the approximation $uu_x \approx u_i^n \frac{u^n_{i+1}-u^n_{i-1}+u^{n+1}_{i+1}-u^{n+1}_{i-1}}{4\Delta x}$, the result will be
\begin{equation}
\begin{aligned}
    & \frac{\alpha}{4}(u^{n+1}_{i+2} - 2u^{n+1}_{i+1} + 2u^{n+1}_{i-1} - u^{n+1}_{i-2}) + u^{n+1}_i - \frac{3\beta}{8}u^n_i(u^{n+1}_{i+1} - u^{n+1}_{i-1})
    \\ & = u^n_i[1+\frac{3\beta}{8}(u^n_{i+1} - u^n_{i-1})] + \frac{\alpha}{4}[u^n_{i-2} - u^n_{i+2} + 2u^n_{i+1} - 2u^n_{i-1}].
\end{aligned}
\end{equation}

However, if I use the approximation $uu_x \approx u_i^{n+1} \frac{u^n_{i+1}-u^n_{i-1}+u^{n+1}_{i+1}-u^{n+1}_{i-1}}{4\Delta x}$, the result will become
\begin{equation}
\begin{aligned}
    & \frac{\alpha}{4}(u^{n+1}_{i+2} - 2u^{n+1}_{i+1} + 2u^{n+1}_{i+1} - u^{n+1}_{i-2}) + \\ & u^{n+1}_i[1-\frac{3\beta}{8}(u^{n+1}_{i+1} - u^{n+1}_{i-1} + u^{n}_{i-1} - u^{n}_{i+1})] 
    \\ & = u^n_i + \frac{\alpha}{4}(u^n_{i-2} - u^n_{i+2} + 2u^n_{i+1} - 2u^n_{i-1})
\end{aligned}
\end{equation}

After gaining the results, we can use matrices to solve these linear equations. Let me first focus on the case of approximation $uu_x \approx u_i^n \frac{u^n_{i+1}-u^n_{i-1}+u^{n+1}_{i+1}-u^{n+1}_{i-1}}{4\Delta x}$. In order to solve the above equations, I view the equations as linear and re-organize them as follows:
\begin{equation}
    \begin{aligned}
        & (\frac{-\alpha}{4})u^{n+1}_{i-2} + (\frac{\alpha}{2}+\frac{3\beta}{8}u^n_i)u^{n+1}_{i-1} + u^{n+1}_i + (-\frac{\alpha}{2} - \frac{3\beta}{8}u^n_i)u^{n+1}_{i+1} + \frac{\alpha}{4}u^{n+1}_{i+2} \\ & = \frac{\alpha}{4}u^n_{i-2} + (-\frac{\alpha}{2} - \frac{3\beta}{8}u^n_i) u^n_{i-1} + u^n_i + (\frac{\alpha}{2}+\frac{3\beta}{8}u^n_i)u^n_{i+1} - \frac{\alpha}{4} u^n_{i+2}.
    \end{aligned}
\end{equation}\label{equ: results_kdv}

Based on the above results, we can re-formulate the problem with the matrix
\begin{equation}
    AU^{n+1} = BU^n,
\end{equation}
where the $A,B,U^n, U^{n+1}$ are all matrices. To be more specific, I write down the details here
\begin{equation}
    \begin{aligned}
        & U^{n+1} = [u^{n+1}_{i-2}\ \ u^{n+1}_{i-1}\ \ u^{n+1}_{i}\ \ u^{n+1}_{i+1}\ \ u^{n+1}_{i+2}]^T, \\ & U^{n} = [u^{n}_{i-2}\ \ u^{n}_{i-1}\ \ u^{n}_{i}\ \ u^{n}_{i+1}\ \ u^{n}_{i+2}]^T,
    \end{aligned}
\end{equation}
Let $\gamma = \frac{\alpha}{2}+\frac{3\beta}{8}u^{n}_i$ for convenience. And the matrices $A,B$ can be represented as
\begin{equation}
\begin{aligned}
& A = \left(
\begin{array}{ccccc}
    1 & -\gamma & \frac{\alpha}{4} & 0 & 0 \\
    \gamma & 1 & -\gamma & \frac{\alpha}{4} & 0 \\
    -\frac{\alpha}{4} & \gamma & 1 & -\gamma & \frac{\alpha}{4} \\
    0 & -\frac{\alpha}{4} & \gamma & 1 & -\gamma \\
    0 & 0 & -\frac{\alpha}{4} & \gamma & 1
\end{array}
\right),
\ 
\\ & B = \left(
\begin{array}{ccccc}
    1 & \gamma & -\frac{\alpha}{4} & 0 & 0 \\
    -\gamma & 1 & \gamma & -\frac{\alpha}{4} & 0 \\
    \frac{\alpha}{4} & -\gamma & 1 & \gamma & -\frac{\alpha}{4} \\
    0 & \frac{\alpha}{4} & -\gamma & 1 & \gamma \\
    0 & 0 & \frac{\alpha}{4} & -\gamma & 1
\end{array}
\right).
\end{aligned}
\end{equation}
Then, I can use formulate the equation as $Ax = b$, where $x = U^{n+1}, b = BU^n$. This is equation can be easily solved in programming. For example, I can use the Python numpy linear algebra library to solve it (i.e., x = np.linalg.solve(A,b)).

\subsection{Stability}
After solving the KdV equation, I need to see whether the Crank-Nicolson scheme is unconditionally stable for this problem. In order to validate its stability of this scheme, I use the Von Neumann stability analysis.

Firstly, I assume $u_i^n = \lambda^n e^{ikj\Delta x}$, where $\lambda$ is a complex number and $k_i$ is the wave number at grid point $i$. Considering this transformation in the solution I calculated above, I can gain the discrete equation. And divide two parts by $\lambda^ne^{ikj\Delta x}$, I can gain the initial results as follows:

Set $$r = \frac{\alpha}{2}+\frac{3\beta}{8} \lambda^n e^{ikj\Delta x},$$

\begin{equation}
    \begin{aligned}
        & \frac{\alpha}{4}\lambda (e^{i2k\Delta x} - e^{-i2k\Delta x}) + \lambda + (-1)r\lambda(e^{ik\Delta x} - e^{-ik\Delta x}) \\ & = \frac{\alpha}{4}\lambda (e^{-i2k\Delta x} - e^{i2k\Delta x}) + 1 + r(e^{ik\Delta x} - e^{-ik\Delta x}),
    \end{aligned}
\end{equation}

Directly solve this equation is hard for me, but I can use the Euler's formula here:
\begin{equation}
    sinx = \frac{e^{ix} - e^{-ix}}{2i}, \ cosx = \frac{e^ix+e^{-ix}}{2}.
\end{equation}
Based on this, I can further transform the equation to the following form:
\begin{equation}
    \lambda = \frac{-0.5\alpha i\ sin(2k\Delta x)+ 2ri\ sin(k\Delta x) + 1}{0.5 \alpha i\ sin(2k \Delta x) - 2ri \ sin(k\Delta x)+1},
\end{equation}

Set $$Q = 0.5 \alpha i\ sin(2k \Delta x) - 2ri \ sin(k\Delta x)+1,$$
\begin{equation}
    \lambda = \frac{2}{Q} - 1.
\end{equation}
Now I just need to check the value of $Q$ to see whether $|\lambda|$ is unconditionally stable. Rewrite the formula of $Q$ as
\begin{equation}
    Q = 2ri \ [sin(k\delta x) - \frac{\alpha}{4r}sin(2k\delta x)] + 1,
\end{equation}
It is cristal clear that $r>0$, $\alpha>0$, $\frac{\alpha}{4r} < 1$, and therefore, I can obtain the result $Q>1$. For any $Q$, I have $|\frac{2}{Q} - 1|<1$. Therefore, I derive that 
\begin{equation}
    |\lambda|<1.
\end{equation}
It is safe to conclude that the Crank-Nicolson scheme in this KdV problem is unconditionally stable.

\subsection{Consistency}

In order to analyze the consistency of the Crank-Nicolson scheme for problem, I need to determine if this method converges to the continuous solution as the discretization parameters (time step $\Delta t$ and spatial step $\Delta x$) approach zero.

Let us begin with the Taylor series of function $u$ in $t$ and $x$, I have the following formulas:

\begin{equation}
\begin{aligned}
    & u^{n+1}_i = u^n_i + ku^n_{t,i}+\frac{k^2}{2}u^n_{tt,i}+\mathcal{O}(k^3),
    \\ & u^{n+1}_{i\pm 1} = u^n_i \pm ku^n_{t,i} + \frac{k^2}{2}u^n_{tt,i}+\mathcal{O}(k^3),
    \\ & 
    u^{n+1}_{i\pm 2} = u^n_i \pm 2ku^n_{t,i} + 2^2\frac{k^2}{2}u^n_{tt,i}+\mathcal{O}(k^3);
    \\ & 
    u^n_{i\pm 1} = u^n_i \pm hu_x + \frac{h^2}{2}u_{xx} - \frac{h^3}{6}u_{xxx} + \mathcal{O}(h^4),
    \\ & 
    u^n_{i\pm 1} = u^n_i \pm 2hu_x + 2^2\frac{h^2}{2}u_{xx} \pm 2^3\frac{h^3}{6}u_{xxx} + \mathcal{O}(h^4),
\end{aligned}
\end{equation}
where $k = \Delta t$ and $h = \Delta x$. Then, replace the relevant terms in Equation \ref{equ: results_kdv} with the above Taylor expansion series, I can obtain:

\begin{equation}
\begin{aligned}
    & u^n_i + (k-\frac{3k\beta u^n_i}{4})u^n_{t,i} + \frac{k^2}{2}u^n_{tt,i} + \mathcal{O}(k^3) \\ & 
    = u^n_i + \frac{3\beta h u^n_i}{4}u_x + (\frac{-\alpha h^3}{2} + \frac{3\beta u^n_i}{8})u_xxx + \mathcal{O}(h^4).
    \end{aligned}
\end{equation}
Then, I define the error function $Error(x,t)$, which use the original KdV equation to minus the above equation. And the error function can be written as:
\begin{equation}
\begin{aligned}
    Error(x,t) = & (1-k+\frac{3\beta u^n_i}{4})u_t + (1 - \frac{\alpha h^3}{2} + \frac{h^3 \beta u^n_i}{8})u_{xxx} \\ & + (\frac{3h\beta}{4} - 1.5) uu_x - \frac{k^2}{2}u_{tt} + \mathcal{O}(k^3+h^4).
\end{aligned}
\end{equation}
Notice that in the previous definitions, I have $k = \Delta t, h = \Delta x, \alpha = \frac{\Delta t}{\Delta x^3}, \beta = \frac{\Delta t}{\Delta x}$. And then, I use these terms in the above $Error(x,t)$ functions:
\begin{equation}
\begin{aligned}
    Error(x,t) = &  (\frac{4 - 4\Delta t + 3\Delta t u^n_i/\Delta x}{4})u_t + (\frac{8-4\Delta t + \Delta t\Delta x^2 u^n_i}{8})u_{xxx} \\ & + (\frac{3\Delta t - 6}{4})uu_x - \frac{\Delta t^2}{2}u_{tt} + \mathcal{O}(k^3 + h^4).
    \end{aligned}
\end{equation}
Based on the above result, when $\Delta t \rightarrow 0, \Delta x \rightarrow 0$, the error function $Error(x,t) \rightarrow 0$. Therefore, the Crank-Nicolson scheme is consistent for the KdV equation here, and I expect it to converge to the exact solution as I refine the discretization parameters.

\subsection{Convergence}
Convergence means that the numerical solution obtained using the scheme converges to the exact solution of the KdV equation as the time step and the spatial discretization are refined.

In order to demonstrate the convergence of the proposed scheme, I need to demonstrate that this finite difference method can be written as a two-level scheme. It follows that I need to demonstrate the matrix $A$ is unconditionally invertible.

Given that I have re-formulated the problem with the matrix in the above sections:
\begin{equation}
    AU^{n+1} = BU^n,
\end{equation}
where the $A,B,U^n, U^{n+1}$ are all matrices (vectors). Specifically,
\begin{equation}
    \begin{aligned}
        & U^{n+1} = [u^{n+1}_{i-2}\ \ u^{n+1}_{i-1}\ \ u^{n+1}_{i}\ \ u^{n+1}_{i+1}\ \ u^{n+1}_{i+2}]^T, \\ & U^{n} = [u^{n}_{i-2}\ \ u^{n}_{i-1}\ \ u^{n}_{i}\ \ u^{n}_{i+1}\ \ u^{n}_{i+2}]^T,
    \end{aligned}
\end{equation}
Let $\gamma = \frac{\alpha}{2}+\frac{3\beta}{8}u^{n}_i$, $\frac{\Delta t}{\Delta x^3} = \alpha,\  \frac{\Delta t}{\Delta x} = \beta$ for convenience. And the matrice $A$ can be represented as
\begin{equation}
\begin{aligned}
& A = \left(
\begin{array}{ccccc}
    1 & -\gamma & \frac{\alpha}{4} & 0 & 0 \\
    \gamma & 1 & -\gamma & \frac{\alpha}{4} & 0 \\
    -\frac{\alpha}{4} & \gamma & 1 & -\gamma & \frac{\alpha}{4} \\
    0 & -\frac{\alpha}{4} & \gamma & 1 & -\gamma \\
    0 & 0 & -\frac{\alpha}{4} & \gamma & 1
\end{array}
\right)_{n \times n},
\end{aligned}
\end{equation}
I can observe that the matrix $A$ is both square ($n \times n$) and symmetric. Now I will try to compute the eigenvalues of this matrix. First, let me formulate this problem as
\begin{equation}
    det(A - \lambda I) = 0,
\end{equation}
where the matrix $I$ is the identical matrix and $\lambda$s are the eigenvalues of the matrix $A$. Optimize the equation above further, I can obtain the following results:
\begin{equation}
\begin{aligned}
& A - \lambda I = \left(
\begin{array}{ccccc}
    1 - \lambda & -\gamma & \frac{\alpha}{4} & 0 & 0 \\
    \gamma & 1 - \lambda & -\gamma & \frac{\alpha}{4} & 0 \\
    -\frac{\alpha}{4} & \gamma & 1 - \lambda & -\gamma & \frac{\alpha}{4} \\
    0 & -\frac{\alpha}{4} & \gamma & 1 - \lambda & -\gamma \\
    0 & 0 & -\frac{\alpha}{4} & \gamma & 1 - \lambda
\end{array}
\right)_{n \times n}.
\end{aligned}
\end{equation}

However, this matrix has a shape of $n \times n$ and has complex elements inside, which makes it hard to directly apply the characteristic equation to obtain the eigenvalues.

The first idea is to use the Jacobi Method, which is an iterative algorithm for determining the solutions of a diagonally dominant system of linear equations. However, it is unsure whether the matrix $A$ is diagonally dominant. Therefore, I consider using the \textit{power iteration method} to gain the eigenvalues.

The iterative power method is stated as follows:
\begin{itemize}
    \item It starts with an initialized vector $b_0$, which can be an approximation to the dominant eigenvector or a random vector;
    \item Update the vector $b_{k+1}$ with $b_k$, and normalize the equation as
    \begin{equation}
        b_{k+1} = \frac{A b_k}{||A b_k||}, \ 0\leq k \leq Max\ Iterations.
    \end{equation}
\end{itemize}
After iterative computations on vector $b$, we can get the eigenvalues of matrix $A$ in $b_{k+1}$ with very small errors.

In our experiments, $\beta = 1,\ \alpha = 1000, \ \gamma = 500+\frac{3}{8}u^n_i, \ \frac{\alpha}{4} = 250$. Therefore, the matrix $A$ will be
\begin{equation}
\begin{aligned}
\left(
\begin{array}{ccccc}
    1 & -500-\frac{3}{8}u^n_i & 250 & 0 & 0 \\
    500+\frac{3}{8}u^n_i & 1 & -500-\frac{3}{8}u^n_i & 250 & 0 \\
    -250 & 500+\frac{3}{8}u^n_i & 1 & -500-\frac{3}{8}u^n_i & 250 \\
    0 & -250 & 500+\frac{3}{8}u^n_i & 1 & -500-\frac{3}{8}u^n_i \\
    0 & 0 & -250 & 500+\frac{3}{8}u^n_i & 1
\end{array}
\right)_{n \times n}.
\end{aligned}
\end{equation}

And applying the iterative power method to this matrix, I can obtain the eigenvalue results. Now I attached the python programming output of iterative power method:

\lstset{
    numbers=left, 
    numberstyle= \tiny, 
    keywordstyle= \color{ blue!70},
    commentstyle= \color{red!50!green!50!blue!50}, 
    frame=shadowbox,
    rulesepcolor= \color{ red!20!green!20!blue!20} ,
    escapeinside=``,
    xleftmargin=0.1em, aboveskip=0.5em,
    framexleftmargin=2em
}
\noindent

\begin{lstlisting}
[0.00470123 0.00354254 0.01583532 ... 0.01614174 
0.00512636 0.00396702]

eigenvalue.shape = (4000,1)
eigenvalues are all positive.
\end{lstlisting}

Note that this result is not just the output of a $5 \times 5$ matrix. In fact, the matrix $A$ has a shape of $4000 \times 4000$. This means the common function that is used to compute eigenvalues such as the $numpy.eig()$ will fail.

Therefore, I have proved that the eigenvalues of our matrix $A$ are unconditionally positive, which indicates that this matrix $A$ is invertible given the fact that the matrix $A$ is both square and symmetric. This also indicates that the proposed scheme can converge on the KdV equation.

\subsection{Considering Another Case}
Let me consider another linearization case in the original KdV equation. That is, linearize the non-linear term $uu_x$ as $uu_x = u^{n+1}_i u_x$. Now putting this term into the KdV equation, I have
\begin{equation}
\begin{aligned}
    & \frac{u^{n+1}_i - u^n_i}{\Delta t} + \frac{1}{4\Delta x^3}[(u^n_{i+2} - 2u^n_{i+1} + 2u^n_{i-1} - u^n_{i-2})  \\ &+ (u^{n+1}_{i+2} - 2u^{n+1}_{i+1} + 2u^{n+1}_{i-1} - 2u^{n+1}_{i-2})] \\ & 
    - 1.5u^{n+1}_i \frac{1}{4\Delta x}(u^{n}_{i+1} - u^{n}_{i-1} + u^{n+1}_{i+1} - u^{n+1}_{i-1}) = 0.
    \end{aligned}
\end{equation}
For convenience, let
\begin{equation}
    \frac{\Delta t}{\Delta x^3} = \alpha,\ \frac{\Delta t}{\Delta x} = \beta.
\end{equation}
Therefore, I can obtain the following results
\begin{equation}
\begin{aligned}
    & \frac{\alpha}{4}u^{n+1}_{i+2} + (\frac{-\alpha}{2} - \frac{3\beta u^{n+1}_i}{8}) u^{n+1}_{i+1} + (1+\frac{3\beta}{8}u^n_{i+1}-\frac{3\beta}{8}u^n_{i-1})u^{n+1}_i \\ & + (\frac{\alpha}{2}+\frac{3\beta}{8}u^{n+1}_i)u^{n+1}_{i-1} - \frac{\alpha}{4} u^{n+1}_{i-2} \\ & = \frac{-\alpha}{4} u^n_{i+2} + \frac{\alpha}{2}u^n_{i+1} + u^n_i - \frac{\alpha}{2}u^n_{i-1} + \frac{\alpha}{4} u^n_{i-2}.
    \end{aligned}
\end{equation}
To better solve the equation above, let me rewrite the equation into the matrix format:
\begin{equation}
    AU^{n+1} = BU^n,
\end{equation}
where the $A,B,U^n, U^{n+1}$ are all matrices. To be more specific, I write down the details here
\begin{equation}
    \begin{aligned}
        & U^{n+1} = [u^{n+1}_{i-2}\ \ u^{n+1}_{i-1}\ \ u^{n+1}_{i}\ \ u^{n+1}_{i+1}\ \ u^{n+1}_{i+2}]^T, \\ & U^{n} = [u^{n}_{i-2}\ \ u^{n}_{i-1}\ \ u^{n}_{i}\ \ u^{n}_{i+1}\ \ u^{n}_{i+2}]^T,
    \end{aligned}
\end{equation}
Let $\zeta = \frac{\alpha}{2}+\frac{3\beta}{8}u^{n+1}_i,\ \eta = 1 + \frac{3\beta u^n_{i+1}}{8} - \frac{3\beta u^n_{i-1}}{8}$. And the matrices $A, B$ can be represented as
\begin{equation}
\begin{aligned}
& A = \left(
\begin{array}{ccccc}
    \eta & -\zeta & \frac{\alpha}{4} & 0 & 0 \\
    \zeta & \eta & -\zeta & \frac{\alpha}{4} & 0 \\
    -\frac{\alpha}{4} & \zeta & \eta & -\zeta & \frac{\alpha}{4} \\
    0 & -\frac{\alpha}{4} & \zeta & \eta & -\zeta \\
    0 & 0 & -\frac{\alpha}{4} & \zeta & \eta
\end{array}
\right),
\\ & B = \left(
\begin{array}{ccccc}
    1 & \frac{\alpha}{2} & -\frac{\alpha}{4} & 0 & 0 \\
    -\frac{\alpha}{2} & 1 & \frac{\alpha}{2} & -\frac{\alpha}{4} & 0 \\
    \frac{\alpha}{4} & -\frac{\alpha}{2} & 1 & \frac{\alpha}{2} & -\frac{\alpha}{4} \\
    0 & \frac{\alpha}{4} & -\frac{\alpha}{2} & 1 & \frac{\alpha}{2} \\
    0 & 0 & \frac{\alpha}{4} & -\frac{\alpha}{2} & 1
\end{array}
\right).
\end{aligned}
\end{equation}
Therefore, the equation can also be written as $Ax = b$, where $x = U^{n+1}$ and $b = BU^n$. However, I notice that the matrix $A$ is not linear for vector $U^{n+1}$. Therefore, it is not appropriate to use matrix to solve the equaiton.

\section{Experiments}
\begin{figure*}[t]
\centering
\subfigure[Initial Data]{
\includegraphics[width=0.27\linewidth]{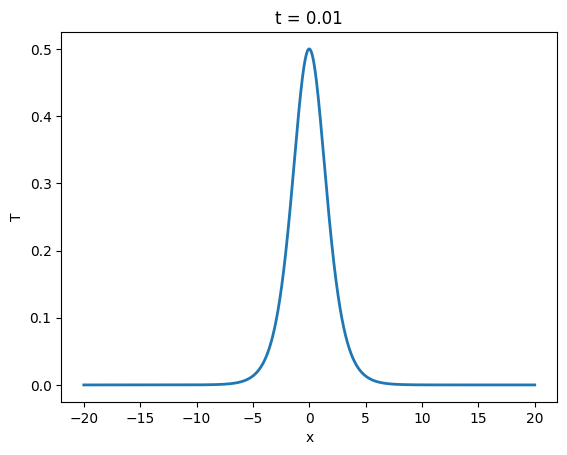}\label{fig:initial}
}
\subfigure[t = 1.01]{
\includegraphics[width=0.27\linewidth]{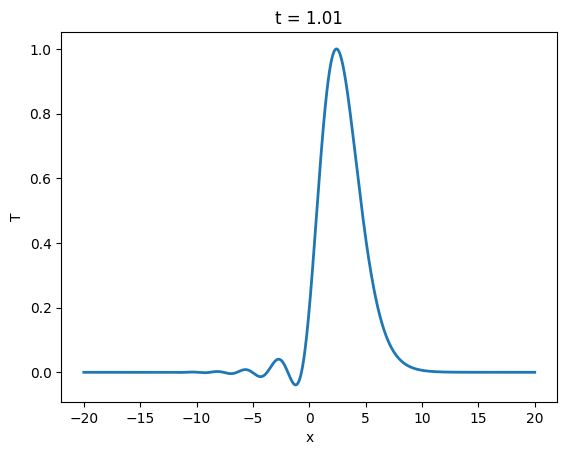}\label{fig:t1}
}
\subfigure[t = 2.01]{
\includegraphics[width=0.27\linewidth]{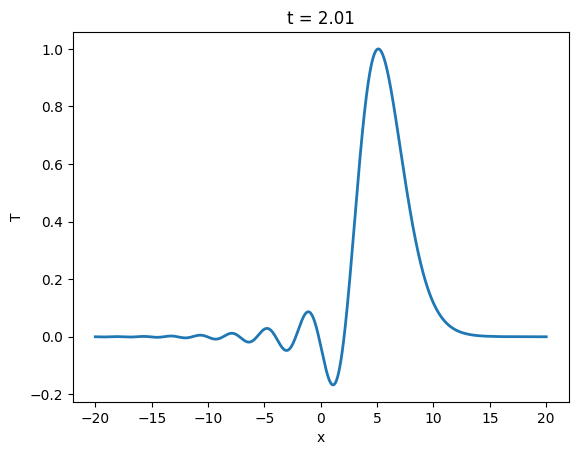}\label{fig:t2}
}
\subfigure[t = 3.01]{
\includegraphics[width=0.27\linewidth]{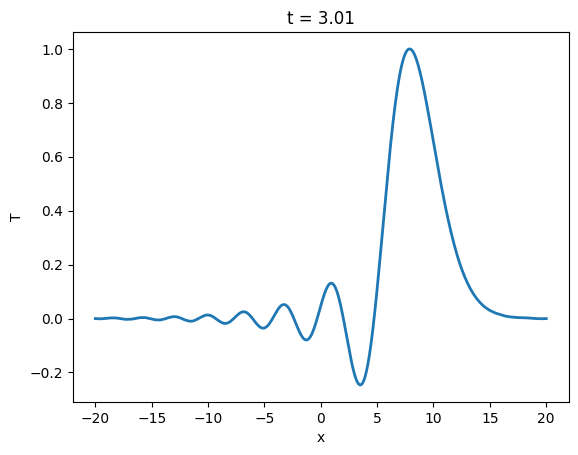}\label{fig:t3}
}
\subfigure[t = 4.01]{
\includegraphics[width=0.27\linewidth]{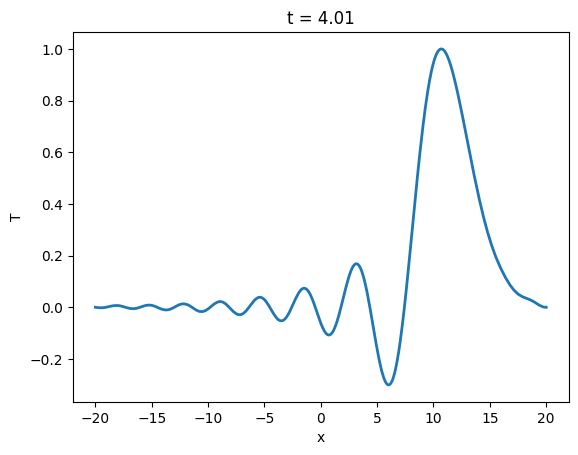}\label{fig:t4}
}
\subfigure[t = 5.01]{
\includegraphics[width=0.27\linewidth]{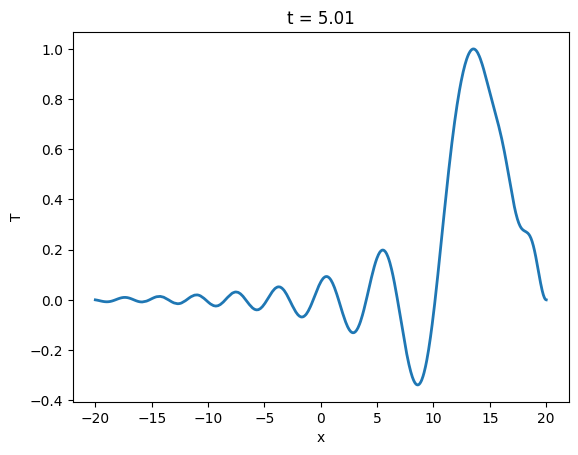}\label{fig:t5}
}
\subfigure[t = 6.01]{
\includegraphics[width=0.27\linewidth]{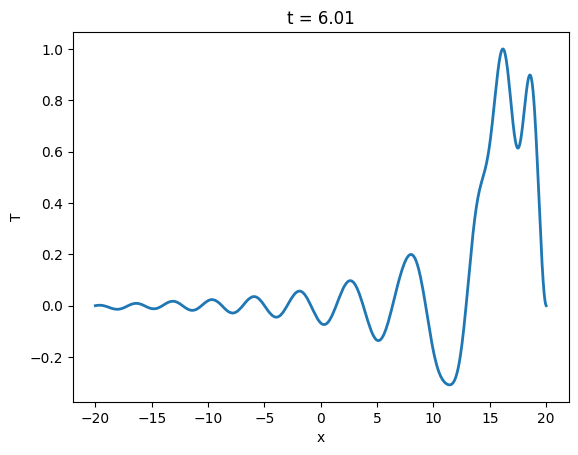}\label{fig:t6}
}
\subfigure[t = 7.01]{
\includegraphics[width=0.27\linewidth]{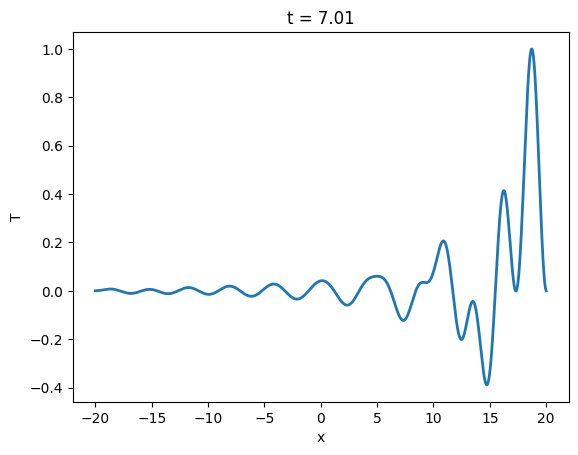}\label{fig:t7}
}
\subfigure[t = 8.01]{
\includegraphics[width=0.27\linewidth]{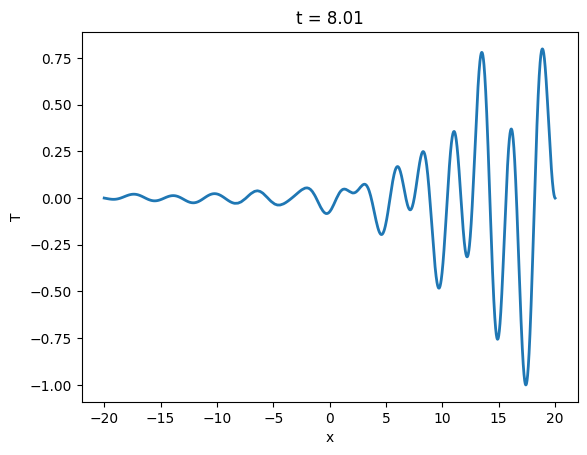}\label{fig:t8}
}\caption{Experimental results of iteratively solving the KdV equation.}\label{fig:experiment}
\end{figure*}

Based on the above theoretical analyses, I implement the python programming method to iteratively solve the KdV equation. And the results are summarized in Figure \ref{fig:experiment}. I plot the 2D figure of the solution of the KdV equation. Note that Figure \ref{fig:initial} represents the initial status $u(x,0)$ of the equation. From the figures, I can see that when time increases, the curves start to become ``flat" and the center of the curve gradually shifts to the right side of number 0. This is a reasonable finding since the solution of the KdV equation has been proved to be $c\times sech^2(\frac{\sqrt(c)}{2})(x-ct)$, which has a shifting center $x = ct$. This experiment also demonstrates the effectiveness of the proposed numerical method.

\section{Conclusion}
In this project report, I have explored the Crank-Nikolson Scheme to solve the KdV equation. By using this implicit finite difference method, I can solve the KdV equation in a much fast way. Also, I show how to derive the scheme on the KdV equation problem, do the stability analysis to prove its unconditional stability, do the consistency analysis and finally combine all these analyses to prove the convergence of the scheme. However, since this is my first time implementing the finite difference method on the partial derivative equations, there may exist some minor mistakes in this report. Helpful discussions or suggestions are warmly welcomed!

\section{Acknowledgement}
Thanks Prof. Bjorn Birnir for his constructive feedback on my report. There are two main differences between the midterm report and this final report:
\begin{itemize}
    \item Complete the ``convergence" section. In this report, I demonstrate the matrix $A$ is unconditionally invertible, which means that this method can be written as a two-level scheme.
    \item Fix the programming code and re-plot the experimental figures (Figure \ref{fig:experiment}). In the midterm report, my figures have many oscillations and I fix this problem by re-setting the parameters. The new code is also attached.
\end{itemize}

\section{Reference}
Also, I would like to list the references that help me with the project. Since the majority of them are online resources, I will directly cite their links.

\begin{itemize}
    \item Stanford University Course (\url{https://web.stanford.edu/class/cme306/Discussion/Discussion1.pdf}).
    \item University of Minnesota (\url{https://www-users.cse.umn.edu/~arnold//papers/stability.pdf}).
    \item "Numerical Solution of Kortweg-de Vries Equation" (\url{https://www.scirp.org/journal/paperinformation.aspx?paperid=99788}).
    \item University of Toronto (\url{https://www.atmosp.physics.utoronto.ca/people/tshaw/tashaw_PHY410_proj.pdf}).
    \item Our textbook of this course (\url{https://epubs.siam.org/doi/book/10.1137/1.9780898717938}).
    \item Youtobe tutorial (\url{https://www.youtube.com/watch?v=y0C3ew3tk2A}).
    \item Online course (\url{https://www.youtube.com/watch?v=qFJGMBDfFMY&list=PLkZjai-2Jcxn35XnijUtqqEg0Wi5Sn8ab}).
\end{itemize}

\section{Appendix}
In this section, I will present my Python coding for the project.

\lstset{
    numbers=left, 
    numberstyle= \tiny, 
    keywordstyle= \color{ blue!70},
    commentstyle= \color{red!50!green!50!blue!50}, 
    frame=shadowbox,
    rulesepcolor= \color{ red!20!green!20!blue!20} ,
    escapeinside=``,
    xleftmargin=0.5em, aboveskip=0.5em,
    framexleftmargin=2em
}
\noindent

\begin{lstlisting}
import numpy as np
import matplotlib.pyplot as plt

# Set up the grid
# L = 100
# T = 10
# N = 1000
# dx = L/N # 0.1
# dt = 0.01
# x = np.arange(-L/2, L/2+dx, dx)
# t = np.arange(0, T+dt, dt)

L = 100
T = 10
N = 10000
dx = L/N
dt = 0.01
x = np.arange(-L/5, L/5+dx, dx)
t = np.arange(0, T+dt, dt)

# Setup the initial and boundary conditionsss
nx = len(x)#1000
nt = len(t)#501

boundaryconditions = [0,0]
initialconditions = 0.5/(np.cosh(x/2)**2)

T = np.zeros((nx,nt))
T[0,:] = boundaryconditions[0]
T[-1,:] = boundaryconditions[1]

T[:,0] = initialconditions

alpha = dt/(dx**3) # 10
beta = 0.375*dt/dx # 0.0375


def make_AB(nx,alpha,beta,un):
     A = np.diag(-0.25*alpha*np.ones(nx-3), -2)+
     np.diag((0.5*alpha+beta*un)*np.ones(nx-2), -1)+ 
     np.diag(np.ones(nx-1), 0)+ 
     np.diag((-0.5*alpha-beta*un)*np.ones(nx-2), 1)+ 
     np.diag(0.25*alpha*np.ones(nx-3), 2)
     
     B = np.diag(0.25*np.ones(nx-3), -2)- 
     np.diag((0.5*alpha + beta*un)*np.ones(nx-2), -1)+ 
     np.diag(np.ones(nx-1), 0)+ 
     np.diag((0.5*alpha+beta*un)*np.ones(nx-2), 1)+ 
     np.diag(-0.25*alpha*np.ones(nx-3), 2)
     
     return A,B

def normalization(data):
    _range = np.max(abs(data))
    return data / _range
    
for j in range(0,nt-1):
    b = T[0:-1,j].copy()
    idx = int((b.shape[0]-1)/2)
    un = b[idx]
    
    A,B = make_AB(nx,alpha,beta,un)

    b = np.dot(B,b) # (999,)
    solution = np.linalg.solve(A,b)
    solution = normalization(solution)
    # update T
    T[1:-1,j+1] = solution[0:-1]

    if j%100 == 0:
        plt.plot(x,T[:,j],linewidth=2)
        plt.ylabel("T")
        plt.xlabel("x")
        plt.title("t = %2.2f"%(dt*(j+1)))
        plt.show()
        plt.clf()
\end{lstlisting}

\end{document}